\documentclass[12pt]{article}

\usepackage{amsmath}
\usepackage{amssymb}
\usepackage{amsthm}
\usepackage{xypic}
\usepackage{url}
\usepackage{dsfont}

\begin{document}
\baselineskip=16pt
\textheight=9.3in
\parindent=0pt 
\def\sk {\hskip .5cm}
\def\skv {\vskip .12cm}
\def\cos {\mbox{cos}}
\def\sin {\mbox{sin}}
\def\tan {\mbox{tan}}
\def\intl{\int\limits}
\def\lm{\lim\limits}
\newcommand{\frc}{\displaystyle\frac}
\def\xbf{{\mathbf x}}
\def\fbf{{\mathbf f}}
\def\gbf{{\mathbf g}}

\def\dbA{{\mathbb A}}
\def\dbB{{\mathbb B}}
\def\dbC{{\mathbb C}}
\def\dbD{{\mathbb D}}
\def\dbE{{\mathbb E}}
\def\dbF{{\mathbb F}}
\def\dbG{{\mathbb G}}
\def\dbH{{\mathbb H}}
\def\dbI{{\mathbb I}}
\def\dbJ{{\mathbb J}}
\def\dbK{{\mathbb K}}
\def\dbL{{\mathbb L}}
\def\dbM{{\mathbb M}}
\def\dbN{{\mathbb N}}
\def\dbO{{\mathbb O}}
\def\dbP{{\mathbb P}}
\def\dbQ{{\mathbb Q}}
\def\dbR{{\mathbb R}}
\def\dbS{{\mathbb S}}
\def\dbT{{\mathbb T}}
\def\dbU{{\mathbb U}}
\def\dbV{{\mathbb V}}
\def\dbW{{\mathbb W}}
\def\dbX{{\mathbb X}}
\def\dbY{{\mathbb Y}}
\def\dbZ{{\mathbb Z}}

\def\la{{\langle}}
\def\ra{{\rangle}}
\def\phi{{\varphi}}

\large\bf\centerline{Fractional Calculus: A Commutative Method on Real Analytic Functions}\rm\normalsize
\vskip .5cm
\centerline{Matthew Parker}
\vskip .75cm

\skv
\large{\bf Abstract}\normalsize
\vskip .1cm
\qquad The traditional first approach to fractional calculus is via the Riemann-Liouville differintegral $_{a}D_{x}^{k}$ [1].  The intent of this paper will be to create a space $K$, pair of maps $g: C^{\omega}(\mathbb{R}) \to K$ and $g': K \to C^{\omega}(\mathbb{R}$), and operator $D^{k}: K \to K$ such that the operator $D^{k}$ commutes with itself, the map $g$ embeds $C^{\omega}(\mathbb{R}$) isomorphically into $K$, and the following diagram commutes;

\begin{displaymath}
    \xymatrix{ C^{\omega}(\mathbb{R}) \ar[d]_{_{a}D_{x}^{k}} \ar[r]^{g} & K \ar[d]^{D^{k}} \\
               C^{\omega}(\mathbb{R}) & K \ar[l]^{g'} }
\end{displaymath} \\*

\qquad This implies the following diagram commutes, for analytic $f$ such that $_{a}D_{x}^{j}f$ = 0 (i.e, if $f = \sum_{i \in I}b_{i}$(x-$a)^{i}$, where \{$b_{i}\} \subset \mathbb{R}$, and $I \subseteq \{j-1, ... , j-\lfloor j \rfloor$\});
\begin{displaymath}
    \xymatrix{f \ar@/_3pc/[dd]_{_{a}D_{x}^{j+k}} \ar[d]^{_{a}D_{x}^{j}} \ar[r]^{g} & g(f) \ar[d]^{D^{j}} \\
               0 & \ar[l]^{g'} D^{j}g(f) \ar[d]^{D^{k}} \\
               _{a}D_{x}^{j+k}f &\ar[l]^{g'} D^{k}D^{j}g(f)  }
\end{displaymath}
\vskip 1cm

\skv
\large{\bf Convention}\normalsize
\vskip .1cm
\qquad Henceforth, unless otherwise noted we assume all functions are real-analytic, thus equal to their Taylor series on some interval of $\mathbb{R}$. When a base point for a Taylor series is not given, we assume it converges on $\mathbb{R}$ or the function has been analytically continued. We let $C^{\omega}(\mathbb{R}$) denote the space of real analytic functions.
\vskip 1cm

\skv
\large{\bf The Space $\mathbb{Z}_{\omega}(a$)}\normalsize
\vskip .1cm
\qquad From basic real analysis, for any $f \in C^{\omega}(\mathbb{R}$), the Taylor series of $f$ (henceforth denoted $T(f$)), equal to $f$ on some open interval in $\mathbb{R}$, is defined by $T(f) = \sum_{i=0}^{\infty}\frac{f^{(i)}(a)}{i!}(x-a)^{i}$ for some $a \in \mathbb{R}$. Note the collection \{$f^{(i)}(a) : i \in \mathbb{Z}_{\geq0}$\} together with the point $a$ uniquely define $f$ within $C^{\omega}(\mathbb{R}$). Then, for fixed $a \in \mathbb{R}$, there is a natural bijection between functions $\sigma: \mathbb{Z}_{\geq0} \to \mathbb{R}$ such that $\sum_{i=0}^{\infty}\frac{\sigma(i)}{i!}(x-a)^{i}$ converges on some interval about $a$, and functions $f \in C^{\omega}(\mathbb{R}$) equal to their Taylor series on some interval about $a$. \\*

$\qquad$ Define $\mathbb{Z}_{\omega}(a$) to be the set of all functions $\sigma : \mathbb{Z} \to \mathbb{R}$ such that $\sum_{i=0}^{\infty}\frac{\sigma(i)}{i!}(x-a)^{i}$ converges on some interval about $a$. When the point $a$ is understood, or not central to the argument but assumed to be fixed, we may omit it and just write $\mathbb{Z}_{\omega}$. Note $\mathbb{Z}_{\omega}$ is non-empty, since the function $f: i \mapsto$ 0 (all $i \in \mathbb{Z}$) is an element of $\mathbb{Z}_{\omega}$. Moreover, $\mathbb{Z}_{\omega}$ is a vector space with identity $\mathds{1}_{\omega}: i \mapsto$ 0; if $\sigma, \rho \in \mathbb{Z}_{\omega}$ and $k \in \mathbb{R}, \sum_{i=0}^{\infty}\frac{(\sigma+\rho)(i)}{i!}(x-a)^{i} = \sum_{i=0}^{\infty}\frac{\sigma(i)}{i!}(x-a)^{i} + \sum_{i=0}^{\infty}\frac{\rho(i)}{i!}(x-a)^{i} \in C^{\omega}(\mathbb{R}$), and $\sum_{i=0}^{\infty}\frac{k\sigma(i)}{i!}(x-a)^{i} = k\sum_{i=0}^{\infty}\frac{\sigma(i)}{i!}(x-a)^{i} \in C^{\omega}(\mathbb{R}$)
\vskip 2cm

\skv
\large{\bf The $R$ Operator}\normalsize
\vskip .1cm
\qquad For $\sigma \in \mathbb{Z}_{\omega}(a$), define the operator $R: \mathbb{Z}_{\omega}(a) \to C^{\omega}(\mathbb{R}$) by $R\sigma = \sum_{i=0}^{\infty}\frac{\sigma(i)}{i!}(x-a)^{i}$. Clearly $R$ is surjective, with kernel \{$\sigma \in \mathbb{Z}_{\omega} : \sigma$(i) = 0, $i \geq$ 0\}. For some pair ($f, a) \in C^{\omega}(\mathbb{R})\times\mathbb{R}$, we also define the operator $R^{-1}: C^{\omega}(\mathbb{R})\times\mathbb{R} \to \mathbb{Z}_{\omega}(a$) by\begin{displaymath}
    R^{-1}(f, a)(i)  = \left\{
     \begin{array}{lr}
       0 & : i\, \textless\, 0\\
       f^{(i)}(a) \, & : i\, \geq\, 0
     \end{array}
   \right.
\end{displaymath} 
This allows the following identities;\\*
\\*
(R1) \,\,$RR^{-1}(f, a) = (f, a$) \\*
(R1') $RR^{-1}f = f$ \\*
(R2) \,\,$R^{-1}R\sigma(i$) =  $\left\{
     \begin{array}{lr}
       \sigma(i) &  \,i \geq 0\\
       0 &  i \;\,\textless\, 0
     \end{array}
   \right.$
$\qquad$ By definition of the maps $R$ and $R^{-1}$, it follows they are both homomorphisms, where $R^{-1}$ is injective with image $\mathbb{Z}_{\omega}$/ker($R$), and $R$ is surjective.
\vskip 4cm

\skv
{\bf The $D$-Operator and $\Gamma$ Function}
\vskip .1cm
\qquad Given $\sigma \in \mathbb{Z}_{\omega}(a$), we define $D^{k}\sigma(i) = \sigma(i+k$) for all $k \in \mathbb{Z}$. From the definition of the operator $D$, we immediately have the identities

(D1) $D^{a}D^{b} = D^{b}D^{a}$ \\*
(D2) $D^{a}D^{b} = D^{a+b}$ \\*
(D3) $D^{a}D^{-a} = D^{-a}D^{a} = D^{0}$ \\*
(D4) $D^{a}(\sigma + \rho) = D^{a}\sigma + D^{a}\rho$ \\*
(D5) $D^{a}(k\sigma) = kD^{a}\sigma$ for all $k \in \mathbb{R}$
\vskip 1cm

\skv
\large{\bf Relating the $D$ Operator to Differentiation}\normalsize
\vskip .1cm
$\qquad$ Induction on the power rule provides the identity $\frac{d^{a}}{dx^{a}}x^{k} = \frac{k!}{(k-a)!}x^{k-a}$, and the relation $n! = \Gamma(n$+1) provides the identity $\frac{d^{a}}{dx^{a}}x^{k} = \frac{\Gamma(k+1)}{\Gamma(k+1-a)}x^{k-a}$. Applying these to Taylor series, we obtain the identity $$T(f) \, = \, \sum_{i=0}^{\infty}\frac{f^{(i)}(a)}{i!}(x-a)^{i} \, = \, \sum_{i=0}^{\infty}\frac{f^{i}(a)}{\Gamma(i+1)}(x-a)^{i}$$ while the power rule allows for the identities $$\frac{d^{j}}{dx^{j}}T(f) \; = \; T(\frac{d^{j}}{dx^{j}}f) \; = \;  \sum_{i=0}^{\infty}\frac{i!}{(i-j)!}f^{(i)}(a)(x-a)^{i-j} \; =  \; \sum_{i=0}^{\infty}\frac{\Gamma(i+1)}{\Gamma(i+1-j)}f^{(i)}(a)(x-a)^{i}$$

$\qquad$ For $f \in C^{\omega}(\mathbb{R}$) and $\sigma \in \mathbb{Z}_{\omega}$ such that $R\sigma = f$, straightforward calculation yields the following identities for all $a \in \mathbb{Z}$; \\*
\\*

(D6) $RD^{a}R^{-1}f = \frac{d^{a}}{dx^{a}}f = f^{(a)}$ \\*
\\*

(D7) $R^{-1}\frac{d^{a}}{dx^{a}}R\sigma(i$) = $\left\{
     \begin{array}{lr}
       \sigma(i+a) & : i \geq -a\\
       0 & : i \,\textless\, -a
     \end{array}
   \right.$
\\*
\\*
(D8) $\frac{d^{a}}{dx^{a}}RD^{-a}\sigma = f$\\*
\\*
$\qquad$ Together, (D1) - (D8), along with (R1') and (R2) will form the core of our arguments for the rest of the paper.

$\qquad$ Finally, we slightly redefine the operator $R$ based on properties of the $\Gamma$ function. By definition, $R\sigma = \sum_{i=0}^{\infty}\frac{\sigma(i)}{\Gamma(i+1)}(x-a)^{i}$. However, for $i \leq 0, \frac{\sigma(i)}{\Gamma(i+1)}$ = 0 so $\frac{\sigma(i)}{\Gamma(i+1)}(x-a)^{i}$ = 0 and $\sum_{i=-\infty}^{\infty}\frac{\sigma(i)}{\Gamma(i+1)}(x-a)^{i} = \sum_{i=0}^{\infty}\frac{\sigma(i)}{\Gamma(i+1)}(x-a)^{i} = R\sigma$, so from this point on we will define 
\vskip 2.5cm
$$R\sigma \; = \; \sum_{i = -\infty}^{\infty}\frac{\sigma(i)}{\Gamma(i+1)}(x-a)^{i}$$ Clearly, properties (R1), (R1'), (R2) and (D1) - (D8) still hold. \\*
\vskip 1cm

\skv
\large{\bf Mapping $C^{\omega}(\mathbb{R}$) To $\mathbb{R}_{\omega}$ And Back}\normalsize
\vskip .1cm
\qquad For a power function $b(x - s)^{\alpha}$, the Riemann-Liouville derivative $_{a}D_{x}^{k}$ is given by [1]$$_{s}D_{x}^{k}b(x-s)^{\alpha} \; = \; \frac{b}{\Gamma(-k)}\int_{s}^{x}(t-s)^{\alpha}(x-t)^{-k-1}dt \; = \; \frac{\Gamma(\alpha+1)}{\Gamma(\alpha+1-k)}b(x-s)^{\alpha-k}$$ $\qquad$ If, and only if, $k \notin \mathbb{R}$ and $\alpha+1-k \in \mathbb{Z}_{\leq0}$, then the numerator of the fraction $\frac{\Gamma(\alpha+1)}{\Gamma(\alpha+1-k)}$ is finite while the denominator goes to $\pm \infty$, and in the limit we see $_{s}D_{x}^{k}f$ = 0. This shows that, when restricted to $C^{\omega}(\mathbb{R}$), ker($_{s}D_{x}^{k}) = \{b(x-s)^{\alpha} : b, \alpha \in \mathbb{R}, \alpha+1-k \in \mathbb{Z}_{\leq 0}$\}. If we wish to preserve the identities (D1) - (D8), (R1') and (R2') when generalizing $D^{k}$ to all real $k$, we must define a new operator with a significantly smaller kernel. Note that (D1) - (D8) and (R1'), (R2) are only consistent if ker($D^{k}$) = \{0\}, the zero function in $\mathbb{Z}_{\omega}$. \\*

$\qquad$ To summarize the situation, then, on the one hand we have the (commutative) $D$ operator on elements of $\mathbb{Z}_{\omega}$ which, when coupled with the $R$ operator, allows identities (D1) - (D8), and on the other we have the Riemann-Liouville derivative, which is commutative for analytic functions when the degrees of differentiation under consideration never sum to a nonpositive integer [2].

$\qquad$ We now create maps $f, f$' and a generalization of the $D$ operator to a space $\mathbb{R}_{\omega}$ (to be defined) such that the following diagram commutes;

\begin{displaymath}
    \xymatrix{C^{\omega}(\mathbb{R}) \ar[r]^f \ar[d]_{_{a}D_{x}^{k}} & \mathbb{R}_{\omega} \ar[d]^{D^{k}}\\
              C^{\omega}(\mathbb{R}) & \mathbb{R}_{\omega} \ar[l]_{f'}}
\end{displaymath}

$\qquad$ It will, however, be more convenient to express the map $f$ as a composition of maps $C^{\omega}(\mathbb{R}) \stackrel{R^{-1}}{\longrightarrow}  \mathbb{Z}_{\omega} \stackrel{\iota}{\longrightarrow} \mathbb{R}_{\omega}$ and extending the domain of the operator $R$ to $\mathbb{R}_{\omega}$ so f = $\iota\circ R^{-1}$ and $f' = R$. Our goal, then, will be to define the maps and spaces which make the following diagram commute, while maintaining analogs of (R1') and (R2);
\vskip 2cm

\begin{displaymath}
    \xymatrix{C^{\omega}(\mathbb{R}) \ar[r]^{R^{-1}} \ar[d]_{_{a}D_{x}^{k}} & \mathbb{Z}_{\omega}\ar[r]^{\iota} & \mathbb{R}_{\omega} \ar[d]^{D^{k}} \\
              C^{\omega} & &\ar[ll]_{R} \mathbb{R}_{\omega} }
\end{displaymath}

$\qquad$ Define $\mathbb{R}_{\omega}(a) = \{\rho: \mathbb{R} \to \mathbb{R} : \sum_{i=-\infty}^{\infty}\frac{\Gamma(i+1-k)}{\Gamma(i+1)}\rho(i)(x-a)^{i-k} \in C^{\omega}(\mathbb{R}) \forall k \in \mathbb{R}$\}. By definition, for any $\rho \in \mathbb{R}_{\omega}, \rho\vert_{_{\mathbb{Z}}} \in \mathbb{Z}_{\omega}$. In fact, for any $\rho = \rho(x) \in \mathbb{R}_{\omega}, \rho(x-k)\vert_{_{\mathbb{Z}}} \in \mathbb{Z}_{\omega}$. Observing $D$ is merely a shift operator on $\mathbb{Z}_{\omega}(\mathbb{R}$), we naturally extend $D$ to $\mathbb{R}_{\omega}$ by setting $D^{k}\rho(i) = \rho$(i-k) for all $k \in \mathbb{R}, \rho \in \mathbb{R}_{\omega}$. By definition of $\mathbb{R}_{\omega}, (D^{k}\rho)\vert_{_{\mathbb{Z}}} \in \mathbb{Z}_{\omega}$ for all $k$. This leads to the natural extension of $R$ to $\mathbb{R}_{\omega}$ by $R\rho = R(\rho\vert_{_{\mathbb{Z}}}) = \sum_{i = -\infty}^{\infty}\frac{\rho(i)}{\Gamma(i+1)}(x-a)^{i}$. \\*

$\qquad$ Properties (D1) - (D5) still hold for $D$ on $\mathbb{R}_{\omega}$, since elements of $\mathbb{R}_{\omega}$, like those of $\mathbb{Z}_{\omega}$, are functions. Thus, we are only left to define the map $\iota$ and verify its properties. \\*

$\qquad$ Let $\sigma \in \mathbb{Z}_{\omega}$, then define $\iota(\sigma$)(z) = lim$_{k \to \infty} (_{a}D_{x}^{z+k}RD^{-k}\sigma$)(a) whenever the limit exists. We then have the following (equivalent) identities;\\*

(I1) $\iota(\sigma)\vert_{_{\mathbb{Z}}} = \sigma$\\*

(I2) ($D^{k}\iota(\sigma))\vert_{_{\mathbb{Z}}} = D^{k}\sigma$ when $k \in \mathbb{Z}$\\*

(I3) $\iota(R^{-1}f)\vert_{_{\mathbb{Z}}} = R^{-1}f$\\*

(I4) $R(\iota(R^{-1}f)\vert_{_{\mathbb{Z}}}) = f$\\*

$\qquad$ Identity (I4) is our analog of (R2), and (R1') follows from properties of $R$ and (I1). Finally, we will show Diagram 2 commutes; that is, $_{a}D_{x}^{k}f = RD^{k}\iota(R^{-1}f$) for all $k \in \mathbb{R}$. Let $f \in C^{\omega}(\mathbb{R}$), and $f_{k}: \mathbb{R} \to \mathbb{R}$ be defined by $f_{k}(z) = (_{a}D_{x}^{z+k}(RR^{-1}f))(a$), then

\begin{align*}
RD^{k}\iota(R^{-1}f) & = R\iota(R^{-1}f - k) \\
 & = R(f_{k}) \\
 & = \sum_{i=0}^{\infty}\frac{1}{\Gamma(i+1)}(_{a}D_{x}^{-i+k}f)(a)(x-a)^{i} \\
 & = T(_{a}D_{x}^{k}f) \\
 & =  _{a}D_{x}^{k}f
\end{align*}
\vskip 3cm

and the diagram commutes. Then (D7) and the following analogs of (D6) and (D8) hold; \\*

(D6') $RD^{k}\iota(R^{-1}f) = _{a}D_{x}^{k}f$\\*

(D8') $_{a}D_{x}^{k}RD^{-k}\iota(R^{-1}f) = f$.\\*
\vskip 1cm

\skv
\large{\bf Conclusion}\normalsize
\vskip .1cm
\qquad  In conclusion, we have created a space $\mathbb{R}_{\omega}$, and a collection of maps and operators $R, R^{-1}, D^{k}$, and $\iota$ such that the operator $D^{k}$ acts exactly the same as the Riemann-Liouville operator as $_{s}D_{x}^{k}$ when applied to an element of $C^{\omega}(\mathbb{R}$) mapped through $\mathbb{R}_{\omega}$, and the operator $D^{k}$ commutes with itself. That is, if $f \in C^{\omega}(\mathbb{R}$) is such that $_{s}D_{x}^{j}f$ = 0 for some $j \in \mathbb{R}$, then for $\sigma = R^{-1}f, \rho = \iota(\sigma$), and all $k \in \mathbb{R}$, we have the following commutative diagram

\begin{displaymath}
    \xymatrix{ f \ar@/_3pc/[dd]_{_{a}D_{x}^{j+k}} \ar[d]^{_{a}D_{x}^{j}} \ar[r]^{R^{-1}} & \sigma \ar[r]^{\iota} & \rho \ar[d]^{D^{j}} \\
               0 & & \ar[ll]^{R} D^{k}\rho \ar[d]^{D^{k}} \\
               _{a}D_{x}^{j+k}f & &\ar[ll]^{R} D^{j}D^{k}\rho  }
\end{displaymath}

which is equivalent to the second diagram in the abstract. This, together with the second diagram in this section - which is equivalent to the first diagram in the abstract - completes the paper.

\end{document}